\documentclass[11pt,reqno]{amsart}
\usepackage{amssymb,verbatim,enumerate,ifthen}
\usepackage[mathscr]{eucal}
\usepackage[utf8]{inputenc}
\usepackage[T1]{fontenc}

\usepackage{dsfont}

\def\N{\mathbb{N}}
\def\R{\mathbb{R}}
\def\Q{\mathbb{Q}}

\def\A{\mathscr{A}}
\def\B{\mathscr{B}}

\def\D{\mathscr{D}}
\def\E{\mathscr{E}}
\def\EE{\mathbb{E}}
\def\F{\mathscr{F}}
\def\FF{\mathbb{F}}
\def\M{\mathscr{M}}
\def\comment#1{}

\def\conv{\mathop{\mbox{\rm conv}}\nolimits}
\def\fix{\mathop{\mbox{\rm Fix}}\nolimits}

\def\ipl{\langle}
\def\ipr{\rangle}

\def\sgn{\mathop{\mbox{\rm sgn}}\nolimits}
\def\one{\mathds{1}}

\newtheorem{theorem}{Theorem}[section]
\newtheorem*{theorem*}{Theorem}
\long\def\Thm#1#2{\ifthenelse{\equal{#1}{*}}{\begin{theorem*}#2\end{theorem*}}
             {\begin{theorem}\label{T#1}#2\end{theorem}}}
\newtheorem{Atheorem}{Theorem}

\def\thm#1{Theorem~\ref{T#1}}

\newtheorem{proposition}[theorem]{Proposition}
\newtheorem*{proposition*}{Proposition}
\long\def\Prp#1#2{\ifthenelse{\equal{#1}{*}}{\begin{proposition*}#2\end{proposition*}}
             {\begin{proposition}\label{P#1}#2\end{proposition}}}
\def\prp#1{Proposition~\ref{P#1}}

\newtheorem{corollary}[theorem]{Corollary}
\newtheorem*{corollary*}{Corollary}
\long\def\Cor#1#2{\ifthenelse{\equal{#1}{*}}{\begin{corollary*}#2\end{corollary*}}
             {\begin{corollary}\label{C#1}#2\end{corollary}}}
\def\cor#1{Corollary~\ref{C#1}}

\newtheorem{lemma}[theorem]{Lemma}
\newtheorem*{lemma*}{Lemma}
\long\def\Lem#1#2{\ifthenelse{\equal{#1}{*}}{\begin{lemma*}#2\end{lemma*}}
             {\begin{lemma}\label{L#1}#2\end{lemma}}}
\def\lem#1{Lemma~\ref{L#1}}

\theoremstyle{definition}
\newtheorem{definition}[theorem]{Definition}
\newtheorem*{definition*}{Definition}
\long\def\Defn#1#2{\ifthenelse{\equal{#1}{*}}{\begin{definition*}\rm #2\end{definition*}}
             {\begin{definition}\label{D#1}\rm #2\end{definition}}}
\def\defn#1{Definition~\ref{D#1}}

\newtheorem{remark}[theorem]{Remark}
\newtheorem*{remark*}{Remark}
\long\def\Rem#1#2{\ifthenelse{\equal{#1}{*}}{\begin{remark*}\rm #2\end{remark*}}
             {\begin{remark}\label{R#1}\rm #2\end{remark}}}

\newtheorem{example}{Example}
\newtheorem*{example*}{Example}
\long\def\Exa#1#2{\ifthenelse{\equal{#1}{*}}{\begin{example*}\rm #2\end{example*}}
             {\begin{example}\label{Ex#1}\rm #2\end{example}}}

\def\eq#1{{\rm(\ref{E#1})}}

\def\Eq#1#2{\ifthenelse{\equal{#1}{*}}
  {\begin{equation*}\begin{aligned}#2\end{aligned}\end{equation*}}
  {\begin{equation}\begin{aligned}\label{E#1}#2\end{aligned}\end{equation}}}

\begin{document}

\date{\today}

\title{Reducible means and reducible inequalities}

\author[T. Kiss]{Tibor Kiss}
\author[Zs. P\'ales]{Zsolt P\'ales}
\address{Institute of Mathematics, University of Debrecen,
H-4032 Debrecen, Egyetem t\'er 1, Hungary}
\email{\{kiss.tibor,pales\}@science.unideb.hu}


\thanks{This research has been supported by the Hungarian
Scientific Research Fund (OTKA) Grant K111651}

\begin{abstract}
It is well-known that if a real valued function acting on a convex set satisfies the $n$-variable 
Jensen inequality, for some natural number $n\geq 2$, then, for all $k\in\{1,\dots, n\}$, it 
fulfills the $k$-variable Jensen inequality as well. In other words, the arithmetic mean and the 
Jensen inequality (as a convexity property) are both reducible. Motivated by 
this phenomenon, we investigate this property concerning more general means and convexity 
notions. We introduce a wide class of means which generalize the well-known means for arbitrary 
linear spaces and enjoy a so-called reducibility property. Finally, we give a sufficient condition for the 
reducibility of the $(M,N)$-convexity property of functions and also for Hölder--Minkowski type inequalities.
\end{abstract}

\maketitle

\section{Introduction}

The notion of Jensen convex functions (introduced by Jensen \cite{Jen05,Jen06} in 1905) plays a central role 
in the theory of convexity and functional inequalities (cf.\ \cite{BecBel61}, \cite{Bul03}, 
\cite{BulMitVas88}, \cite{HarLitPol34}, \cite{Kuc85}, \cite{MitPecFin93}, \cite{RobVar73}). To recall its 
classical definition, let $D$ be a convex subset of a real linear space $X$. Then we say that the 
function $f:D\to\R$ is 
\emph{Jensen convex} if 
\Eq{J2}{
f\Big(\frac{x_1+x_2}{2}\Big)\leq\frac{f(x_1)+f(x_2)}{2},\qquad(x_1,x_2\in D).
}
It is an important property of Jensen convex functions that, for all $n\in\N$, they also satisfy the 
$n$-variable Jensen inequality
\Eq{Jn}{
f\Big(\frac{x_1+\dots+x_n}{n}\Big)\leq\frac{f(x_1)+\dots+f(x_n)}{n},\qquad(x_1,\dots,x_n\in D).
}
The standard proof of this inequality (based upon \eq{J2}) uses a particular induction which is attributed to 
Cauchy: First, using normal induction on $k$, it is proved that \eq{Jn} holds for $n=2k$. Then, 
assuming that \eq{Jn} holds for $n=k$, it is deduced that it is also valid for $n=k-1$. We call 
this the reducibility property of the $n$-variable Jensen inequality \eq{Jn}.

The idea is to replace the two appearance of the arithmetic mean in \eq{Jn} by arbitrary 
means $M:D^n\to D$ and $N:I^n\to I$, and to consider functions $f:D\to I$ satisfying
\Eq{MN}{
f\big(M(x_1,\dots,x_n)\big)
\leq N\big(f(x_1),\dots,f(x_n)\big),\qquad(x_1,\dots,x_n\in D).
}
Our main aim is to find and describe general sufficient conditions under which, for 
$k\in\{1,\dots,n\}$, a $k$-variable 
convexity property can be deduced from \eq{MN}. This requires the construction of $k$-variable means which 
are the reductions of $M$ and $N$, respectively. The construction and computation of the 
$k$-variable reductions will be elaborated in the class of deviation means introduced by Daróczy 
\cite{Dar71b,Dar72b} (which includes Hölder and Gini means \cite{Gin38}, quasi-arithmetic means 
\cite{HarLitPol34}, Matkowski means \cite{Mat10b} and Bajraktarević means \cite{Baj58,Baj63}), and 
also in the class of generalized deviation means that will be introduced in this paper to provide a 
broad class of means for the vector valued setting. We also demonstrate how generalized deviation 
means can be derived as solutions of convex minimum problems. Finally, we consider and 
establish the reducibility property of Hölder--Minkowski type inequalities under natural 
assumptions.

\section{Terminology and notations}

We adopt the standard notations $\N$, $\Q$, and $\R$ for the sets of natural, rational and real 
numbers, respectively, furthermore $\R_+$ denotes the set of positive real numbers, that is
\Eq{*}{
\R_+:=\,]0,+\infty[\,:=\{t\in\R\mid t>0\}.
}
Given a natural number $n\in\N$, we shall frequently use the notation $\N_n$ defined as
\Eq{*}{
\N_n:=\{1,\dots,n\}:=[1,n]\cap\N.
}
\comment{For a given set of natural numbers $H\subseteq \N$, we define the \emph{characteristic function}
$\one_H:\N\to\{0,1\}$ of $H$ by
\Eq{*}{
\one_H(i):=
\begin{cases}
1 & \text{if }i\in H,\\[2mm]
0 & \text{if }i\in \N\setminus H.
\end{cases}
}}
For an arbitrary nonempty set $S$ and $n\in\N$, we also identify the elements of the Cartesian product 
$S^n$ with the set of all functions mapping $\N_n$ to $S$, that is, with the set 
$S^{\N_n}:=\{x:\N_n\to S\}$. Furthermore, for $x\in S^n$ and $i\in\N_n$, we simply denote $x(i)$ by 
$x_i$.

Finally, we introduce a notation which will be applied throughout this paper. 
Let $n\in\N$, $k\in\N_n$, let $\chi:\N_k\to\N_n$ be an injective function and $S$ be a set.
For $x\in S^k$ and $y\in S$, the symbol $(x|\chi)(y)$ denotes the element of $S^n$ defined by
\Eq{*}{
  (x|\chi)(y)_i
  :=\begin{cases}
  y & \mbox{if } i\in\N_n\setminus\chi(\N_k),\\[1mm]
  x_j & \mbox{if } i\in\chi(\N_k) \mbox{ and } i=\chi(j).
  \end{cases}
}

\section{Reducible means}

In the sequel, let $X$ be a linear space over $\R$ and $D\subseteq X$ be a nonempty convex set. For 
a given $H\subseteq X$, the set $\conv(H)\subseteq X$ denotes the \emph{convex hull} of $H$, namely 
the smallest convex subset of $X$ which contains $H$. It is easy to see, that $\conv(H)$ is the set 
of all vectors in $X$ which can be written as a convex combination of finitely many elements of 
$H$, thus, obviously, $\conv(D)\subseteq D$.

\Defn{mean}{
Let $n\in\N$. We say that an $n$-variable function 
$M:D^n\to X$ is a \emph{mean} on $D$ if
\Eq{*}{
M(x)\in\conv (x(\N_n)),
\qquad\mbox{that is,}\qquad
M(x_1,\dots,x_n)\in\conv\{x_1,\dots,x_n\},
\qquad(x\in D^n).
}
The mean $M$ will be called \emph{strict} if, for all $x\in D^n$, the vector
$M(x)$ belongs to the \emph{relative interior} of $\conv\{x_1,\dots,x_n\}$, that is, $M(x)$ can be 
written as convex combination of $x_1,\dots,x_n$ with positive coefficients. We say that $M$ is 
\emph{symmetric} if, for all bijection $p:\N_n\to\N_n$, we have
\Eq{*}{
M(x\circ p)=M(x),
\qquad\mbox{that is,}\qquad
M(x_{p_1},\dots,x_{p_n})=M(x_1,\dots,x_n),
\qquad(x\in D^n).
}}

It immediately follows from this definition and from the convexity of $D$ that, for a mean 
$M:D^n\to X$, we always have $M(D^n)\subseteq D$, and that $M$ is \emph{reflexive}, which means 
that $M(u,\dots,u)=u$ holds for all $u\in D$. 

The most important example of an $n$-variable mean is the \emph{arithmetic mean} $\A:X^n\to X$ defined by
\Eq{A}{
\A(x)=\A(x_1,\dots,x_n):=\frac{x_1+\cdots+x_n}{n}.
}
More generally, if $\omega:X\to\R_+^n$, then the \emph{functionally weighted arithmetic mean} 
$\A^\omega:X^n\to X$ is defined by
\Eq{*}{
 \A^\omega(x)
 :=\A\Big(\!\!
\begin{array}{c}
x_1\quad,\dots,\quad x_n \\ 
\omega_1(x_1),\dots,\omega_n(x_n) 
\end{array}
\!\!\Big)
:=\frac{\omega_1(x_1)x_1+\cdots+\omega_n(x_n)x_n}{\omega_1(x_1)+\cdots+\omega_n(x_n)}.
}
If the function $\omega$ is constant on $X$, then we simply speak about a \emph{weighted arithmetic 
mean}. One can easily see that $\A^\omega$ is a strict mean and it is a symmetric mean if 
$\omega_1=\cdots=\omega_n$. 

More general means will be constructed in terms of deviations and families of convex functions in the next 
section.

In what follows, we define the notions of continuity and reduction of mean $M:D^n\to X$ with respect to a 
given injective map $\chi:\N_k\to\N_n$.

\Defn{cont}{
Let $M:D^n\to X$ be a mean. We say that $M$ is \emph{$\chi$-continuous} if, for any $x\in D^k$, the 
mapping $m_{x,M}:\conv (x(\N_k))\to X$ defined as
\Eq{mxy}{
m_{x,M}(y):=M\big((x|\chi)(y)\big).
}
is continuous on $\conv (x(\N_k))$.
}

\Defn{red}{
Let $M:D^n\to X$ be a mean. We say that $M$ is \emph{$\chi$-reducible} if there exists a mean 
$K:D^k\to X$ such that, for all $x\in D^k$, the vector $y=K(x)$ is a solution of the equation
\Eq{e1}{
M\big((x|\chi)(y)\big)=y.
}
The mean $K$ will be called a \emph{$\chi$-reduction of} $M$. If for all $x\in D^k$, the equation 
\eq{e1} has a unique solution $y\in\conv(x(N_k))$, that is, if $K$ is uniquely determined, then
we say that $M$ is a \emph{uniquely $\chi$-reducible mean}, furthermore, the mean $K$ will be called 
the \emph{$\chi$-reduction of $M$} and will be denoted by $M_\chi$.
}

The next theorem is about the existence of $\chi$-reductions.

\Thm{1}{
If the mean $M:D^n\to X$ is $\chi$-continuous, then it is also $\chi$-reducible.}

\begin{proof}
Let $x\in D^k$ be arbitrarily fixed and define the function $m_{x,M}:\conv (x(\N_k))\to X$ by \eq{mxy}.
Obviously, the target set of $m_{x,M}$ is $\conv (x(\N_k))$, and, because of the $\chi$-continuity of the 
mean $M$, the function $m_{x,M}$ is continuous on the compact convex set $\conv (x(\N_k))$. Thus, due to the 
Brouwer Fixed Point Theorem, the fixed point set
\Eq{*}{
\fix(m_{x,M}):=\{y\in\conv (x(\N_k))\mid m_{x,M}(y)=y\}
}
is not empty. Finally, define $K(x)$ to be any element of the nonempty set $\fix(m_{x,M})$. Then, for all 
$x\in D^k$, the vector $y=K(x)$ will be a solution of \eq{e1}, hence $K$ is a $\chi$-reduction of $M$.
\end{proof}

For the setting of unique $\chi$-reducibility, we shall need the following useful lemma.

\Lem{4}{
Let $I\subseteq\R$ be an interval, $n\in\N$, $k\in\N_n$, and $\chi:\N_k\to\N_n$
be an injective function. Assume that the $\chi$-continuous mean $M:I^n\to\R$ 
is uniquely $\chi$-reducible. Then, for all $x\in I^k$ and for all $y\in J_x:=[\min(x),\max(x)]$, 
we have
\Eq{mx2}{
\sgn\big(m_{x,M}(y)-y\big)=\sgn(M_{\chi}(x)-y),
}
where $m_{x,M}:\conv (x(\N_k))\to\R$ is defined by \eq{mxy}.}

\begin{proof}
Let $x\in I^k$ be arbitrarily fixed. If $\min(x)=\max(x)$, then the statement is obvious,
thus we may assume that $\min(x)<\max(x)$. For the sake of brevity, define 
\Eq{*}{
\mu_{x,M}(y):=m_{x,M}(y)-y
}
for $y\in J_x$. Then, due to the definition of the $\chi$-reduction 
of means,  we have $\mu_{x,M}(y)=0$ for $y\in J_x$ if and only if $y=M_{\chi}(x)$.

First assume that $M_{\chi}(x)$ belongs to the interior of $J_x$. Because of the mean-property
of $M$, obviously, we have $\mu_{x,M}(\max(x))<0<\mu_{x,M}(\min(x))$. Then, because of the
uniqueness of the zero of $\mu_{x,M}$ and of the $\chi$-continuity of $M$ on the interval $J_x$,
it immediately follows that $\mu_{x,M}$ must be strictly positive on the subinterval 
$[\min(x),M_{\chi}(x)[\,$, and it must be strictly negative on the subinterval  $]M_{\chi}(x),\max(x)]$.

On the other hand, if either $M_{\chi}(x)=\min(x)$ or $M_{\chi}(x)=\max(x)$, then a similar argument 
shows that the function $\mu_{x,M}$ is strictly positive on the interval $J_{x}\setminus\{\min(x)\}$ or it is
strictly negative on the entire interval $J_{x}\setminus\{\max(x)\}$, respectively,
which finishes the proof.
\end{proof}

Note that if an $n$-variable symmetric mean is reducible for some injective function mapping $\N_k$ 
to $\N_n$, then it is also reducible with respect to any injective $\N_n$-valued function defined 
on the set $\N_k$.

The prototypical example for this phenomenon is the arithmetic mean defined in \eq{A}. More 
precisely, the $n$-variable arithmetic mean $A:X^n\to X$ is $\chi$-reducible with respect to any 
injective function $\chi:\N_k\to\N_n$. Indeed, for a fixed $x\in X^k$, the equation \eq{e1} of 
\defn{red} has the form
\Eq{*}{
\frac{x_1+\dots+x_k+(n-k)y}{n}=y.
}
A direct calculation shows that $y:=(x_1+\dots+x_k)/k$ is the only solution of the equation above 
on the entire space $X$. Thus the $k$-reduced $A$ mean of $x\in X^k$ is just its $k$-variable 
arithmetic mean.

For the sake of brevity, we introduce the following notation: If $S$ is an arbitrary nonempty set and 
$u=(u_1,\dots,u_n)\in S^n$ then $u_\chi$ denotes the $k$-tuple $(u_{\chi_1},\dots,u_{\chi_k})\in 
S^k$. Concerning the $\chi$-reduction of a functionally weighted arithmetic mean, we have the 
following result.

\Prp{A}{
Let $\omega:D\to\R_+^n$. Then we have $\A^\omega_\chi=\A^{\omega_\chi}$.
}

\begin{proof}
For the mean $M=\A^\omega$ and for $x\in D^k$, equation \eq{e1} can be rewritten as
\Eq{*}{
\frac{\omega_{\chi_1}(x_1)x_{1}+\dots+\omega_{\chi_k}(x_k)x_k+\big(\sum_{i\,\not\in\,\chi(\N_k)}
\omega_i(y)\big)y}
{\omega_{\chi_1}(x_1)+\dots+\omega_{\chi_k}(x_k)+\sum_{i\,\not\in\,\chi(\N_k)}\omega_i(y)}=y.
}
It immediately follows that the unique solution $y$ of this equation is of the form
\Eq{*}{
  y=\frac{\omega_{\chi_1}(x_1)x_{1}+\dots+\omega_{\chi_k}(x_1)x_k}
     {\omega_{\chi_1}(x_1)+\dots+\omega_{\chi_k}(x_k)}=\A^{\omega_\chi}(x),
}
which proves that $\A^\omega_\chi(x)=\A^{\omega_\chi}(x)$.
\end{proof}

\section{Generalized deviation functions and generalized deviation means}

We recall now the notion of standard deviation function and deviation mean,
which was first introduced and investigated by Zoltán Daróczy in \cite{Dar72b}. This class 
of means has many interesting properties (Aczél and Daróczy \cite{AczDar63c}, Daróczy \cite{Dar72b,Dar71b}, 
Daróczy--Losonczi \cite{DarLos70}, Daróczy--Páles \cite{DarPal82,DarPal83}, Losonczi 
\cite{Los71a,Los71b,Los71c,Los73a}, Páles 
\cite{Pal83a,Pal82a,Pal84a,Pal84b,Pal85a,Pal87d,Pal88a,Pal88e,Pal88d}) and it generalizes the well-known 
classes of means (for instance Hölder means \cite{HarLitPol34}, Gini means \cite{Gin38}, quasi-arithmetic 
means \cite{HarLitPol34} and quasi-arithmetic means with weight function, that is, Bajraktarević 
means \cite{Baj58,Baj63}).

Let $I\subseteq\R$ be an interval. A function $E:I\times I\to\R$ is called a \emph{deviation 
function} (shortly a \emph{deviation}) if the following two properties hold: 
\begin{enumerate}[(E1)]\itemsep=2mm
 \item $E(u,u)=0$ for all $u\in I$ and,
 \item for any fixed $u\in I$, the function $v\mapsto E(u,v)$ is continuous and strictly decreasing 
on $I$.
\end{enumerate}
The class of deviation functions defined on $I$ will be denoted by $\EE(I)$.
These two properties imply that, for a deviation $E\in\EE(I)$, we always have the following 
so called sign-property: 
\Eq{sgn}{
\sgn E(u,v)=\sgn(u-v),
\qquad (u,v\in I).
}
Now, using a finite collection of deviations, we can derive means on the interval $I$. For 
$E\in\EE(I)^n$ and $x\in I^n$, the unique value $y\in I$, satisfying the equation
\Eq{dev}{
E_1(x_1,y)+\dots+E_n(x_n,y)=0,
}
is called the \emph{$E$-deviation mean} or \emph{$E$-Daróczy mean} of $x$, and is denoted by 
$\D^E(x)$. Observe that the notion of the deviation mean is well-defined. Indeed, let $x\in I^n$ be 
arbitrarily fixed, and denote $\alpha:=\min(x)$ and $\beta:=\max(x)$. The continuity and the strict 
decreasingness of the function $\E_E:I\to\R$ defined by
\Eq{EE}{
\E_{E,\,x}(u):=E_1(x_1,u)+\dots+E_n(x_n,u),
}
and the inequalities $\E_{E,\,x}(\alpha)\geq0\geq\E_{E,\,x}(\beta)$ show that there uniquely exists 
$y\in[\alpha, \beta]$ such that \eq{dev} holds.

Finally, we have the following easy-to-prove but useful statement.

\Lem{ldev}{
Let $n\in\N$ and $E=(E_1,\dots,E_n)\in\EE(I)^n$. Then, using the notation \eq{EE}, for all $x\in 
I^n$,
\Eq{*}{
\sgn\circ\,\E_{E,\,x}(u)=
\sgn\big(\D^E(x)-u\big),\qquad(u\in I).
}}

We recall now the most classical examples for deviation functions and the means generated by them. Let 
$f:I\to\R$ be a strictly increasing continuous function and $\omega:I\to\R_+$. Then, the 
two-variable function $E_{f,\omega}:I\times I\to\R$, defined by
\Eq{*}{
  E_{f,\omega}(x,y):=\omega(x)(f(x)-f(y)),
}
is trivially a deviation. In the particular case, when $\omega$ equals the constant $1$ and $f$ is 
the identity function, then $E_{f,\omega}$ is called the \emph{arithmetic deviation}. To generate 
more general deviation means, we consider two settings.

\begin{enumerate}[(i)]\itemsep=2mm
\item\emph{Bajraktarević type means.} Let $f:I\to\R$ be a strictly increasing 
continuous function, $\omega:I\to\R_+^n$, and define $E_i:=E_{f,\omega_i}$ for $i\in\N_n$.
Then, for any $x\in I^n$, the equation \eq{dev} has the following explicite solution:
\Eq{*}{
y=\B^{f,\omega}(x):=f^{-1}\bigg(\frac{\omega_1(x_1)f(x_1)+\cdots+\omega_n(x_n)f(x_n)}
    {\omega_1(x_1)+\cdots+\omega_n(x_n)}\bigg),
    }
which is called the \emph{Bajraktarević mean} generated by $f$ and $\omega$. In fact, this mean was 
introduced by Bajraktarević \cite{Baj58} in the particular case when $\omega_1=\cdots=\omega_n$. If $f(x)=x$ 
on $I$, then the mean so obtained is equal to functionally weighted arithmetic mean $\A^\omega$. If, for some 
$p\neq q$, we have that $f(x)=x^{p-q}$ and $\omega_1(x)=\cdots=\omega_n(x)=x^q$, then the above expression 
yields Gini means.
\item\emph{Generalized quasi-arithmetic means or Matkowski type means.} Let $f_1,\dots,f_n:I\to\R$ 
be strictly increasing continuous functions and let $\omega:I\to\R_+^n$ be the constant function 
$\omega(t):=(1,\dots,1)$. Finally, define $E_i:=E_{f_i,\omega_i}$ for $i\in\N_n$. Then, 
for any $x\in I^n$, the solution of equation \eq{dev} can be directly calculated again and has the 
form
\Eq{*}{
y=\M^{f,\omega}(x)
:=(f_1+\dots+f_n)^{-1}(f_1(x_1)+\dots+f_n(x_n)),
}
which is called the \emph{Matkowski mean} generated by $f$ (cf.\ \cite{Mat10b}, \cite{MatPal15}). In the 
particular case when $f_1=\cdots=f_n$, the above expression simplifies to a so-called quasi-arithmetic mean 
which has a rich theory developed in the book \cite{HarLitPol34}. By taking $f_1(x)=\cdots=f_n(x)=x^p$ for 
some nonzero real $p$, Hölder means (or power means) can also be obtained.
\end{enumerate}

In the rest of the paper, let $X$ be a \textit{Hausdorff topological vector space} over 
$\R$. For an arbitrary nonempty subset $S\subseteq X$, let $S^*$ denote the the space of all 
continuous linear functionals defined on the linear hull of $(S-S)$. In what follows, we shall 
extend the notion of deviation function and deviation mean to convex subsets of linear spaces.

\Defn{gd}{
Let $D\subseteq X$ be a nonempty convex set. We say that a mapping $E:D\times D\to D^*$ is a 
\textit{generalized deviation function} if it satisfies the following two properties:
\begin{enumerate}[(GE1)]\itemsep=2mm
\item\label{d0} $E(u,u)=0$ for all $u\in D$, and
\item\label{cm} for all fixed $u\in D$, the function $v\mapsto -E(u,v)$ is 
continuous and  \textit{strictly monotone} on $D$, that is
\Eq{*}{
(E(u,v)-E(u,w))(v-w)<0,
\qquad (u,v,w\in D\text{ with }v\neq w).}
\end{enumerate}
The class of generalized deviation functions defined on $D$ will be denoted by $\EE(D)$.}
Observe that the properties (\ref{d0}) and (\ref{cm}) imply that, for a generalized deviation $E\in\EE(D)$,
we always have 
\Eq{gsgn}{
E(u,v)(u-v)>0,
\qquad (u,v\in D,\,u\neq v).
}

Now, using a finite collection of generalized deviations, we can define means on the convex set $D$. In 
contrast to the definition of deviation means (that are defined on real intervals), the notion of generalized 
deviation mean will be defined by a system of inequalities.

\Defn{gdm}{
Let $E=(E_1,\dots,E_n)\in\EE(D)^n$. For $x\in D^n$, we say that the vector $y\in\conv (x(\N_n))$ is the 
\textit{generalized $E$-deviation mean} of $x$ if
\Eq{fi}{
(E_1(x_1,y)+\dots+E_n(x_n,y))(x_i-y)\leq 0,
\qquad(i\in\N_n).
}
If $y\in\conv (x(\N_n))$ exists and unique, then it will be denoted by $\D^E(x)$.}

The next theorem states that the notion of generalized $E$-deviation mean is well-defined.

\Thm{wdf}{
Let $n\in\N$ and $E\in\EE(D)^n$. Then, for all $x\in D^n$, there 
uniquely exists $y\in\conv(x(\N_n))$ such that \eq{fi} holds.
}

\begin{proof}
Let $x\in D^n$ be arbitrarily fixed and, for the brevity, denote the compact 
convex set $\conv(x(\N_n))$ by $C_x$ and define the function $\E_{E,x}:D\to D^*$ by \eq{EE}.
Then, by the defining properties of generalized deviations, the function $-\E_{E,x}$ is continuous and 
strictly monotone. Observe, that the real valued mapping $\phi:C_x\times C_x\to\R$, given by 
\Eq{*}{
\phi(u,v):=\E_{E,x}(u)(v-u),
}
is continuous in its first variable, and (in view of the linearity of $\E_{E,x}(u)(\cdot)$ for any fixed 
$u\in C_x$) is affine (convex and concave simultaneously) in its second variable. Thus, due to the Ky Fan 
Minimax Inequality Theorem (cf.\ \cite{Cha05}, \cite{Cla13}, \cite{Zei88}), there exists $y\in C_x$, such that
\Eq{*}{
\sup_{v\in C_x}\E_{E,x}(y)(v-y)
 =\sup_{v\in C_x}\phi(y,v)
 \leq\sup_{w\in C_x}\phi(w,w)
 =\sup_{w\in C_x}\E_{E,x}(w)(w-w)=0.
}
Thus, for every $v\in C_x$, in particular, for every $v\in\{x_1,\dots,x_n\}$, we have
\Eq{*}{
  \E_{E,x}(y)(v-y)\leq 0.
}
This proves the existence of $y\in \conv(x(\N_n))$ satisfying \eq{fi}.

To prove the uniqueness, assume indirectly, that there exist $y\neq z$ in $ \conv(x(\N_n))$ satisfying 
\eq{fi}. Then, for all $i\in\N_n$, we have
\Eq{np+}{
\E_{E,x}(y)(x_i-y)\leq 0
\qquad\text{and}\qquad
\E_{E,x}(z)(x_i-z)\leq 0.
}
The vectors $y,z$ being in $\conv(x(\N_n))$, there exist convex combination coefficients 
$\lambda_1,\dots,\lambda_n\geq0$ with $\lambda_1+\dots+\lambda_n=1$ and 
$\mu_1,\dots,\mu_n\geq0$ with $\mu_1+\dots+\mu_n=1$ such that
\Eq{*}{
y=\lambda_1x_1+\dots+\lambda_nx_n
\qquad\text{and}\qquad
z=\mu_1x_1+\dots+\mu_nx_n.
}
Multiplying the first and second inequalities in \eq{np+} by $\mu_i$ and $\lambda_i$, respectively, and then 
adding up the inequalities so obtained, we get
\Eq{*}{
\E_{E,x}(y)(z-y)\leq 0
\qquad\text{and}\qquad
\E_{E,x}(z)(y-z)\leq 0.
}
The sum of these two inequalities can be written as
\Eq{np*}{
\big(\E_{E,x}(y)-\E_{E,x}(z)\big)(y-z)\geq0.
}
On the other hand, using the strict monotonicity of $(-\E_{E,x})$, we obtain that
\Eq{*}{
  (\E_{E,x}(y)-\E_{E,x}(z))(y-z)<0,
}
which contradicts \eq{np*}. This proves that the vector $y\in \conv(x(\N_n))$, satisfying the inequality 
\eq{fi}, is uniquely determined.
\end{proof}

\Rem{dgd}{It is obvious that if $X:=\R$ and $D\subseteq\R$ is an interval, then $D^*\equiv\R$ and the notion 
of generalized deviation
functions and generalized deviation means reduces to that of deviation functions and deviation means, 
respectively.

To verify the statement about the means, let $n\in\N$, $E\in\EE(D)^n$, and $x\in D^n$ be 
arbitrary, and assume that $\min(x)<\max(x)$. We need to show that the value 
$y\in D$ is the solution of the equation \eq{dev} in $D$ if and 
only if it is the solution of the system of inequalities \eq{fi} in $\conv(x(\N_n))=[\min(x),\max(x)]$.

If the value $y\in D$ is the solution of \eq{dev}, that is, it is the $E$-deviation mean of $x$, 
then, the inequalities $\E_{E,x}(\min(x))\geq0\geq\E_{E,x}(\max(x))$ show that $y\in[\min(x),\max(x)]$ and it 
trivially satisfies the inequalities of \eq{fi}, that is, the vector $y$ is the generalized $E$-deviation mean of $x$.

Conversely, assume that $y\in[\min(x),\max(x)]$ is the generalized $E$-deviation mean of $x$, or 
equivalently, it is the solution of the system \eq{fi}. Then, in particular, we have
\Eq{mm}{
\E_{E,x}(y)\cdot(\min(x)-y)\leq 0
\qquad\text{and}\qquad
\E_{E,x}(y)\cdot(\max(x)-y)\leq 0.
}
If $y$ were one of the endpoints of the interval $[\min(x),\max(x)]$, say $y=\min(x)$, then $y<\max(x)$, 
therefore the second inequality yields that $\E_{E,x}(y)\leq 0$. On the other hand, $y\leq 
x_i$ for all $i\in\N_n$, and, for at least one index $j\in\N_n$, we have that $y<x_j$. Thus, for all 
$i\in\N_n$, the inequalities $E_i(x_i,y)\geq0$ and $E_j(x_j,y)>0$ hold. This implies that $\E_{E,x}(y)>0$. The 
contradiction so obtained shows $y$ is bigger than $\min(x)$. Similarly, $y$ is smaller than $\max(x)$. 
Therefore, the two inequalities in \eq{mm} result that $\E_{E,x}(y)$ is nonnegative and also 
nonpositive. Consequently, we must have $\E_{E,x}(y)=0$, that is, $y$ is the $E$-deviation mean of $x$.
}

\Thm{rgd}{
Let $n\in\N$, $k\in\N_n$ and let $E\in\EE(D)^n$. Then the generalized $E$-deviation mean $\D^E:D^n\to D$ 
is reducible with respect to any injective function
$\chi:\N_k\to\N_n$. Furthermore, the $\chi$-reduction of $\D^{E}$ is uniquely determined, namely
\Eq{*}{
\D_{\chi}^{E}(x)=\D^{E_\chi}(x),
\qquad (x\in D^k).}}

\begin{proof}
Let $x\in D^k$ be arbitrarily fixed and denote $y_0:=\D^{E_\chi}(x)$. The property (\ref{d0}) of generalized 
deviations implies that
\Eq{*}{
E_i\big((x|\chi)(y)_i,y\big)=
\begin{cases}
0 & \text{if }i\in\N_n\setminus\chi(\N_k),\\[2mm]
E_i(x_{j},y) & \text{if } i\in\chi(\N_k)\text{ and }i=\chi(j).
\end{cases}
}
Therefore, 
\Eq{*}{
\E_{E,(x|\chi)(y)}(y)&=E_1\big((x|\chi)(y)_1,y\big)+\dots+E_n\big((x|\chi)(y)_n,y\big)\\[1mm]
&=E_{\chi_1}(x_1,y)+\dots+E_{\chi_k}(x_k,y)=\E_{E_\chi,x}(y).
}

According to \defn{red}, we need to show that $y=y_0$ is the unique solution of the equation 
$\D^E\big((x|\chi)(y)\big)=y$ 
in $\conv(x(\N_k))$, that is, $y=y_0$ is the unique solution of the system of inequalities
\Eq{*}{
\E_{E,(x|\chi)(y)}(y)\big((x|\chi)(y)_i-y\big)=\E_{E_\chi,x}(y)\big((x|\chi)(y)_i-y\big)\leq 0, 
\qquad(i\in\N_n).
}
The inequalities automatically hold when $i\in\N_n\setminus\chi(\N_k)$ (because then $(x|\chi)(y)_i=y$), 
therefore the above system of inequalities is equivalent to
\Eq{nk}{
  \E_{E_\chi,x}(y)(x_i-y)\leq 0, \qquad(i\in\N_k).
}
In view of \thm{wdf}, the system of inequalities in \eq{nk} is uniquely solvable in $\conv(x(\N_k))$ and its 
$y$ solution equals $y_0=\D^{E_\chi}(x)$, which was to be proved.
\end{proof}

In the theorem below, we construct the large class of generalized deviations in terms families of relatively
G\^ateaux differentiable strictly convex functions. As a consequence of such a representation, generalized 
deviation means can be viewed as the unique minimizers of certain strictly convex functions.

Given an arbitrary set $S\subseteq X$, a point $u\in S$ is called a \emph{relative algebraic 
interior point of $S$} if, for all $v\in S$, the set $\{t\in\R\mid tv+(1-t)u\in S\}$ is a right
neighborhood of $0$ in $\R$. The set $S$ is said to be \emph{relatively algebraically open} if every point of $S$ is 
its relative algebraic interior point.

A function $f:S\to\R$ is called \emph{relatively G\^ateaux differentiable} at a relatively algebraically interior point
$u$ of $S$ if there exists a continuous linear functional $f'(u)\in S^*$ such that, for all $v\in S$, 
\Eq{Gd}{
\lim_{t\to0^+}\frac{f(u+t(v-u))-f(u)}{t}=f'(u)(v-u).
}
The notion of G\^ateaux differentiability with respect to a subspace of $X$ (in our case, with respect to the 
linear span of $S-S$), was considered in the paper \cite{PalZei08}.

We need the following auxiliary result, which is the adaptation of some well-known theorems about convex 
functions to our setting (cf.\ the books \cite{RobVar73} and \cite{Zal02}.

\Thm{AR}{
Let $D\subseteq X$ be a convex set and $f:D\to\R$ be a relatively G\^ateaux differentiable function on $D$. 
Then the following statements hold.
\begin{enumerate}[(1)]\itemsep=1mm
\item $D$ is relatively algebraically open and, for every $u\in D$, the relative G\^ateaux derivative $f'(u)$ 
is uniquely determined. 
\item\label{ar2} The function $f$ is convex if and only if 
\Eq{AR1}{
f(v)\geq f(u)+f'(u)(v-u),
\qquad(u,v\in D),
}
and $f$ is strictly convex if and only if this inequality is strict whenever $u\neq v$.
\item\label{ar3} The function $f$ is convex if and only if its G\^ateaux derivative $f'$ is 
monotone, that is,
\Eq{AR2}{
(f'(u)-f'(v))(u-v)\geq0, \qquad (u,v\in D),
}
and $f$ is strictly convex if and only if this inequality is strict whenever $u\neq v$.
\item If $S\subseteq D$ is a nonempty convex set and $f$ attains its minimum at $u\in S$ 
on the set $S$, then
\Eq{AR0}{
f'(u)(v-u)\geq0,\qquad (v\in S),
}
Conversely, if $f$ is convex and \eq{AR0} holds for some $u\in S$, then $f$ attains its minimum at $u$ on the 
set $S$.
\end{enumerate}}

\begin{proof}
Let $u\in D$ be arbitrarily fixed. Then, because of the convexity of $D$, for all $v\in D$, we have 
$[0,1]\subseteq\{t\in\R\mid tv+(1-t)u\in S\}$, which shows that $u$ is a relative 
algebraic interior point of $D$. Assume that $f'(u)$ is not uniquely determined, that is, there exists 
$\varphi,\psi\in D^*$ such that, for all $v\in D$,
\Eq{}{
\lim_{t\to0^+}\frac{f(u+t(v-u))-f(u)}{t}=\varphi(v-u)=\psi(v-u).
}
Then, $(\varphi-\psi)(v-u)=0$ for all $v\in D$. Now, let $h\in D-D$ be arbitrary. Then there exist $v,w\in D$ such that 
$h=v-w$, hence
\Eq{*}{
  (\varphi-\psi)(h)=(\varphi-\psi)(v-u)-(\varphi-\psi)(w-u)=0.
}
Therefore, $\varphi-\psi$ vanishes on the linear span of $D-D$, showing that $\varphi=\psi$.

To prove (\ref{ar2}), assume that $f$ is convex. Then, for all $u,v\in D$, the map $t\mapsto 
\frac{1}{t}(f(u+t(v-u))-f(u))$ is nondecreasing, hence
\Eq{*}{
  f(v)-f(u)=\frac{f(u+1(v-u))-f(u)}{1}\geq \lim_{t\to0}\frac{f(u+t(v-u))-f(u)}{t}=f'(u)(v-u),
}
which gives \eq{AR1}. If $f$ is strictly convex and $u\neq v$, then $t\mapsto 
\frac{1}{t}(f(u+t(v-u))-f(u))$ is 
strictly increasing, which results that \eq{AR1} holds with strict inequality.

For the converse, assume \eq{AR1}, and let $u,v\in D$ and $t\in[0,1]$ be arbitrary. Then, based 
on \eq{AR1}, we get that
\Eq{AR3}{
  f(u)&\geq f(tu+(1-t)v)+f'(tu+(1-t)v)(u-(tu+(1-t)v))\\
 &=f(tu+(1-t)v)+(1-t)f'(tu+(1-t)v)(u-v), 
\\[2mm]
  f(v)&\geq f(tu+(1-t)v)+f'(tu+(1-t)v)(v-(tu+(1-t)v))\\
  &=f(tu+(1-t)v)+tf'(tu+(1-t)v)(v-u). 
}
Multiplying the first inequality by $t$, the second one by $(1-t)$, and adding up the inequalities so 
obtained  side by side, we get
\Eq{*}{
 tf(u)+(1-t)f(v)\geq f(tu+(1-t)v),
}
which proves the convexity of $f$. If \eq{AR1} holds with strict inequality for $u\neq v$ and $x\neq 
y$, then the inequalities in \eq{AR3} are strict for $t\not\in\{0,1\}$, hence we obtain the 
strict convexity of $f$. 

To prove the second assertion, assume again that $f$ is convex. Then \eq{AR1} holds, thus, applying 
this inequality twice, we obtain that 
\Eq{*}{
f(v)\geq f(u)+f'(u)(v-u)\qquad\mbox{and}\qquad
f(u)\geq f(v)+f'(v)(u-v)
}
for all $u,v\in D$. Adding up these inequalities side by side, it results that \eq{AR2} is valid. 
If $f$ is strictly convex and $u\neq v$, then \eq{AR1} is strict, which yields that \eq{AR2} is also 
strict.

Conversely, assume \eq{AR2} and, for $u,v\in D$, define the function $f_{u,v}:[0,1]\to\R$ by
\Eq{*}{
f_{u,v}(t):=f(tu+(1-t)v).
}
Observe that $f_{u,v}$ is differentiable on $[0,1]$, furthermore the derivative 
$f'_{u,v}:=\frac{d}{dt}f_{u,v}$ is nondecreasing. Indeed, a short calculation shows that
\Eq{*}{
f'_{u,v}(t)=\lim_{\tau\to t}\frac{f_{u,v}(\tau)-f_{u,v}(t)}{\tau-t}=f'(tu+(1-t)v)(u-v).
}
Now let $t,s\in[0,1]$ such that $t\neq s$. Then, due to \eq{AR2}, we have
\Eq{*}{
0\leq(t-s)(f'(tu+(1-t)v)-f'(su+(1-s)v))(u-v)
=(t-s)(f'_{u,v}(t)-f'_{u,v}(s)),
}
which implies that $f'_{u,v}$ is nondecreasing. We obtained that $f_{u,v}$ is convex for any fixed 
$u,v\in D$.

Finally, let $u,v\in D$ and $t\in[0,1]$ be arbitrarily fixed. Then we have the following 
calculation:
\Eq{*}{
f(tu+(1-t)v)=f_{u,v}(t)&=f_{u,v}(t\cdot 1+(1-t)\cdot 0)\\
&\leq tf_{u,v}(1)+(1-t)f_{u,v}(0)=tf(u)+(1-t)f(v),
}
consequently $f$ is convex.

For the third statement, let $S\subseteq D$ be a nonempty convex set and assume that $f$ attains its minimum 
on $S$ at the point $u\in S$. Then, for all $t\in[0,1]$ and $v\in S$, we have that $f(u+t(v-u))\geq f(u)$. 
Hence, in view of formula \eq{Gd}, we get that $f'(u)(v-u)\geq0$ for all $v\in S$. 

Now assume that $f$ is convex and, for some $u\in S$, \eq{AR0} holds. Then, applying \eq{AR1} for $u,v\in S$, 
we get
\Eq{*}{
  f(v)\geq f(u)+f'(u)(v-u)\geq f(u).
}
This proves that $f$ attains its minimum on $S$ at the point $u\in S$.
\end{proof}

To formulate the next theorem, let $\FF(D)$ denote the class of functions $F:D\times D\to\R$ 
with the following property:
\begin{itemize}\itemsep=2mm
\item[(F)] for any fixed $u\in D$, the function $F_u:=F(u,\cdot)$ is relatively G\^ateaux 
differentiable and strictly convex on $D$, furthermore $F_u'(u)=0$.
\end{itemize}

\Thm{ch}{
Assume that $D\subseteq X$ is a convex set and let $F\in\FF(D)$.
Then the function $E_F:D\times D\to D^*$, defined by 
\Eq{ph3}{
 E_F(u,v)=-F_u'(v),
}
is a generalized deviation. Furthermore, if $n\in\N$, $F\in\FF(D)^n$ and 
$E_F=(E_{F_1},\dots,E_{F_n})$, then, for $x\in D^n$, the equality $y=\D^{E_F}(x)$ holds if and only 
if $y$ is the unique minimizer over $\conv(x(\N_n))$ of the function $\F_{F,x}:D\to\R$ defined by
\Eq{ffx}{
  \F_{F,x}(v):=F_1(x_1,v)+\cdots+F_n(x_n,v).
}
Conversely, if $X$ is the real line and $D$ is an open interval, then, for all deviations $E\in\EE(D)$, there 
exists a function $F\in\FF(D)$ such that, for all $u\in D$,
\Eq{fue}{
F_u'(v)=-E(u,v),\qquad(v\in D)
}
is satisfied.}

\begin{proof}
First let $F\in\FF(D)$ and define the function $E_F:D\times D\to D^*$ as in \eq{ph3}. We show that 
$E_F$ is a generalized deviation. It only suffices to verify the strict monotonicity of $-E_F$ in its second 
variable. Let $u,v,w\in D$ such that $v\neq w$. According to the property (F) of $F$, the function 
$F_u$ is strictly convex on its domain, or equivalently, based on \thm{AR}, we have that
\Eq{*}{
0<(F_u'(v)-F_u'(w))(v-w)
=-(E_F(u,v)-E_F(u,w))(v-w).
}
Consequently, the function $-E_F(u,\cdot)$ is strictly monotone on $D$. 

Now let $n\in\N$, $F\in\FF(D)^n$, $E_F=(E_{F_1},\dots,E_{F_n})$ and let $x\in D^n$ be arbitrarily 
fixed. The function $\F_{F,x}:D\to\R$, defined in \eq{ffx}, is continuous on the convex, compact set 
$\conv(x(\N_n))$, thus there exists a point $y\in\conv(x(\N_n))$, which minimizes $\F_{F,x}$ on 
the set $\conv(x(\N_n))$. Moreover, because of the strict convexity of $\F_{F,x}$, the minimizer 
$y$ is unique. Thus, based on the last statement of \thm{AR}, for all $v\in\conv(x(\N_n))$, we have
\Eq{*}{
0\leq\F_{F,x}'(y)(v-y)=-(E_{F_1}(x_1,y)+\dots+E_{F_n}(x_n,y))(v-y).
}
In particular, this inequality holds also for all $v\in\{x_1,\dots,x_n\}$. Because of the 
uniqueness of the generalized $E_F$-deviation mean of $x$ (cf. \thm{wdf}), we must have 
$y=\D^{E_F}(x)$.

Conversely, if 
\Eq{*}{
(E_{F_1}(x_1,y)+\dots+E_{F_n}(x_n,y))(v-y)\leq0
}
for all $v\in\{x_1,\dots,x_n\}$, then this inequality is also valid for all $v\in\conv(x(\N_n))$.
Hence, for all $v\in\conv(x(\N_n))$,
\Eq{*}{
 \F_{F,x}'(y)(v-y)\geq0.
}
In view of the reversed implication in the last statement of \thm{AR}, this implies that $y$
is the minimizer of the function $\F_{F,x}$ over the set $\conv(x(\N_n))$.

Let finally $X:=\R$ and $D\subseteq\R$ be an interval, furthermore let $E\in\EE(D)$ be a deviation 
and define the function $F:D\times D\to\R$ by the formula
\Eq{int}{
F(u,v):=-\int\limits_u^vE(u,t)\,dt,\qquad (u,v\in D).
}
For all $u\in D$, the function $t\mapsto E(u,t)$ is continuous on $D$, thus, due to the Fundamental Theorem 
of Calculus, $F_u$ is continuously differentiable on $D$, and \eq{fue}
holds. The strict decreasingness of $E$ in its second variable implies that $F_u'$ is a strictly monotone 
and hence $F_u$ is strictly convex. Obviously we also have that $F_u'(u)=-E(u,u)=0$ for all $u\in D$.
\comment{Finally, we show that $F$ satisfies the property (\ref{ph1}). 
By \eq{int}, we have that $F(u,u)=0$ for all $u\in D$. Let 
$u,v\in D$ such that $u\neq v$. If $u<v$, then, due to the sign-property \eq{sgn} of
deviations, the function $t\mapsto E(u,t)$ is negative on the interval $]u,v[\,$, thus, according 
again to \eq{int}, we get $F(u,v)>0$. 
On the other hand, if $u>v$, then $t\mapsto E(u,t)$ is positive on $]v,u[\,$, hence
\Eq{*}{
F(u,v)
=-\int\limits_u^v E(u,t)\,dt
=\int\limits_v^u E(u,t)\,dt>0.
}
We obtained that $F(u,v)\geq0$ for all $u,v\in D$ and $F(u,v)=0$ if and only if $u=v$, that is, 
the property (\ref{ph1}) is satisfied.}
\end{proof}

The following result offers the construction of families of strictly convex functions in terms of two 
single variable functions. We recall that the unit ball of a normed space $(X,\|\cdot\|)$ is called 
strictly convex if $\|x\|=\|y\|=1$ and $x\neq y$ implies that $\|tx+(1-t)y\|<1$ for all $t\in]0,1[$.
(Observe that the strict convexity of the unit ball does not imply that the norm is a strictly convex 
function, moreover, by the positive homogeneity, any norm cannot be strictly convex.)

\Prp{1}{Let $(X,\|\cdot\|)$ be a normed space, assume that the unit ball is strictly convex and 
the norm is G\^ateaux differentiable on $X\setminus\{0\}$. Let further $D\subseteq X$ be a convex set and
$\omega:D\to\R_+$. Then the function $F:D\times D\to\R$, defined by
\Eq{F}{
F(u,v):=\omega(u)\|v-u\|^2,\qquad(u,v\in D),
}
satisfies property (F).}

\begin{proof} Let $u\in D$ be fixed. To show that $F_u:=F(u,\cdot)$ is strictly convex, let $v,w\in D$ 
with $v\neq w$ and $t\in\,]0,1[\,$. We distinguish two cases. 

First assume that the vectors $v-u$ and $w-u$ are not parallel (that is, there is no $t\in[0,1]$ such that 
$t(v-u)=(1-t)(w-u)$). Then non of them is zero and $x:=\frac{v-u}{\|v-u\|}$ and $y:=\frac{w-u}{\|w-u\|}$ are 
distinct unit vectors. Therefore, by the strict convexity of the unit ball, we have that 
$\|sx+(1-s)y\|<1$ for all $s\in\,]0,1[\,$. Now, by also using the convexity of the square function, we get
\Eq{*}{
  F_u&(tv+(1-t)w)=\omega(u)\|tv+(1-t)w-u\|^2=\omega(u)\|t(v-u)+(1-t)(w-u)\|^2 \\
    &=\omega(u)\big(t\|v-u\|+(1-t)\|w-u\|\big)^2
    \bigg\|\frac{t\|v-u\|}{t\|v-u\|+(1-t)\|w-u\|}x+\frac{(1-t)\|w-u\|}{t\|v-u\|+(1-t)\|w-u\|}y\bigg\|^2 \\
    &<\omega(u)\big(t\|v-u\|+(1-t)\|w-u\|\big)^2 \\
    &\leq \omega(u)\big(t\|v-u\|^2+(1-t)\|w-u\|^2\big) 
    = tF_u(v)+(1-t)F_u(w).
}
Secondly, assume that $v-u$ and $w-u$ are parallel vectors. Then, the relation $v\neq w$ 
implies that $\|v-u\|\neq\|w-u\|$. Thus, by the subadditivity and the positive homogeneity of the norm and 
the strict convexity of the square function, we get
\Eq{*}{
  F_u(tv+(1-t)w)&=\omega(u)\|t(v-u)+(1-t)(w-u)\|^2 
    \leq\omega(u)\big(t\|v-u\|+(1-t)\|w-u\|\big)^2 \\
    &<\omega(u)\big(t\|v-u\|^2+(1-t)\|w-u\|^2\big) 
    = tF_u(v)+(1-t)F_u(w).
}
To check the G\^ateaux differentiability, denote $p(x):=\|x\|$ and let $v\in D\setminus\{u\}$ and $h\in X$. 
Then
\Eq{*}{
  F_u'(v)(h)
   &=\lim_{t\to0^+}\frac{F_u(v+th)-F_u(v)}{t}
    =\lim_{t\to0^+}\frac{\omega(u)\|v+th-u\|^2-\omega(u)\|v-u\|^2}{t} \\
  &=\omega(u)\lim_{t\to0^+}\big(\|v+th-u\|+\|v-u\|\big)\frac{p(v-u+th)-p(v-u)}{t} \\
  &=2\omega(u)\|v-u\|p'(v-u)(h).
}
Therefore, for $u\neq v$, we get $F_u'(v)=2\omega(u)\|v-u\|p'(v-u)$.

On the other hand, for $v=u$, we have
\Eq{*}{
  F_u'(u)(h)
   &=\lim_{t\to0^+}\frac{F_u(u+th)-F_u(u)}{t}
    =\lim_{t\to0^+}\frac{\omega(u)\|th\|^2-\omega(u)\|0\|^2}{t} 
    =\omega(u)\lim_{t\to0^+}\frac{t^2\|h\|^2}{t}=0,
}
which proves that $F_u'(u)=0$. This completes the proof of property (F).
\end{proof}

\Exa{*}{
Let $(X,\ipl\cdot,\cdot\ipr)$ be an inner product space over $\R$, $D\subseteq X$ be a 
nonempty convex set, and $\omega:D\to\R_+$. Then, by the previous result, the function $F:D\times D\to\R$, 
defined by \eq{F} belongs to $\FF(D)$, and for all $u,v\in D$, we have
\Eq{EF}{
E_{F}(u,v)(h)=-F_u'(v)(h)=-2\omega(u)\|v-u\|p'(v-u)(h)=2\omega(u)\ipl u-v,h\ipr,\qquad(h\in X).
}
Now we can explicitly compute the generalized deviation mean generated by such generalized deviations. Let 
$n\in\N$, $\omega_1,\dots,\omega_n:D\to\R_+$ and $F_1,\dots,F_n:D\times D\to\R$ be functions,
defined as in \eq{F} using the weight functions $\omega_1,\dots,\omega_n$, respectively, furthermore let 
$E_F:=(E_{F_1},\dots, E_{F_n})$. Then
\Eq{*}{
\D^{E_F}(x)=\frac{\omega_1(x_1)x_1+\dots+\omega_n(x_n)x_n}{\omega_1(x_1)+\dots+\omega_n(x_n)}=\A^\omega(x),
\qquad(x\in D^n).
}
Indeed, for $x\in D^n$ and $h\in X$, with the notation $y:=\A^\omega(x)\in\conv(x(\N_n))$, we have
\Eq{*}{
(E_{F_1}(x_1,y)+\dots+E_{F_n}(x_n,y))(h)
&=2(\omega_1(x_1)\ipl x_1-y,h\ipr+\dots+\omega_n(x_n)\ipl x_n-y,h\ipr)\\[2mm]
&=2\ipl\omega_1(x_1)x_1+\dots+\omega_n(x_n)x_n-(\omega_1(x_1)+\dots+\omega_n(x_n))y,h\ipr=0.
}
In particular, this equality holds also for $h\in\{x_1,\dots,x_n\}-y$, thus we must have $y=\D^{E_F}(x)$.

On the other hand, by \thm{ch}, the vector $y=\A^\omega(x)$ is the unique minimizer of the function
\Eq{*}{
 \F_{F,x}(v):=F_1(x_1,v)+\cdots+F_n(x_n,v)=\omega_1(x_1)\|x_1-v\|^2+\cdots+\omega_n(x_n)\|x_n-v\|^2,
}
that is, $y$ is \emph{the weighted least square approximant} of the elements $x_1,\dots,x_n\in D$.}

\section{Reducible inequalities involving means}

In this section we consider convexity properties, comparison and Hölder--Minkowski type inequalities and 
establish their reducibility.

\Defn{mn}{
Let $D\subseteq X$ be a nonempty convex set, $n\in\N$ and let $M:D^n\to X$ and $N:\R^n\to\R$
be means. We say that a function $f:D\to\R$ is \emph{convex with respect to the pair of means $(M,N)$
on $D$} or that $f$ is \emph{$(M,N)$-convex on $D$} if
\Eq{mnc}{
(f\circ M)(x)\leq N(f\circ x),
\qquad(x\in D^n),
}
that is, if
\Eq{*}{
f\big(M(x_1,\dots,x_n)\big)
\leq N\big(f(x_1),\dots,f(x_n)\big),\qquad(x_1,\dots,x_n\in D).
}}

\Thm{rp2}{
Let $D\subseteq X$ be a nonempty convex set, $I\subseteq\R$ be an interval, $n\in\N$, $k\in\N_n$, and let 
$\chi:\N_k\to\N_n$ be an injective function. Let further $M:D^n\to X$ and $N:I^n\to\R$ be means such that 
$M$ is $\chi$-reducible and $N$ is $\chi$-continuous and uniquely $\chi$-reducible. If a function 
$f:D\to I$ is $(M,N)$-convex, then it is also $(K,N_{\chi})$-convex for all $\chi$-reduction 
$K:D^k\to X$ of the mean $M$.}

\begin{proof}
Let $f:D\to I$ be an $(M,N)$-convex function, $K:D^n\to X$ be any $\chi$-reduction of $M$ 
and let $x\in D^k$ be arbitrarily fixed. Denote $y:=K(x)$. Then, because of the definition of 
$y$, we have $M\big((x|\chi)(y)\big)=y$. Using this, the $(M,N)$-convexity of $f$, and the notation \eq{mxy}, 
we obtain that
\Eq{*}{
f(y)=(f\circ M)\big((x|\chi)(y)\big)\leq N\big(f\big((x|\chi)(y)\big)\big)=m_{f\circ x,N}(f(y)),
}
which is equivalent to the inequality
\Eq{*}{
0\leq m_{f\circ x,N}(f(y))-f(y).
}
Due to the $\chi$-continuity and to the unique $\chi$-reducibility of $N$, using \lem{4}, it 
immediately follows that $f(y)\leq N_{\chi}(f\circ x)$ holds, that is
\Eq{*}{
(f\circ K)(x)\leq N_{\chi}(f\circ x).
}
Consequently, $f$ is $(K,N_{\chi})$-convex on its domain.
\end{proof}

The subsequent corollaries immediately follow from the theorem above, from \prp{A} and from 
\thm{rgd}.

\Cor{1K}{Let $D\subseteq X$ be a nonempty convex set, $I\subseteq\R$ be an interval and $n\in\N$. 
Let further $\omega:D\to\R_+^n$ and $E:I\times I\to\R^n$ such that $E_i$ is a deviation for all 
$i\in\N_n$. If a function $f:D\to I$ satisfies the $n$-variable inequality
\Eq{*}{
f\big(\A^\omega(x_1,\dots,x_n)\big)
\leq\D^E\big(f(x_1),\dots,f(x_n)\big),\qquad(x_1,\dots,x_n\in D),
}
then, for all $k\in\N_n$ and for all injective function $\chi:\N_k\to\N_n$, it also satisfies
the $k$-variable inequality
\Eq{*}{
f\big(\A^{\omega_{\chi}}(x_1,\dots,x_k)\big)
\leq\D^{E_{\chi}}\big(f(x_1),\dots,f(x_k)\big),\qquad (x_1,\dots,x_k\in D).
}}

\Cor{2K}{Let $D\subseteq X$ be a nonempty convex set, $I\subseteq\R$ be an interval and $n\in\N$. 
Let further $G:D\times D\to (D^*)^n$ and $E:I\times I\to\R^n$ such that $G_i$ is a generalized 
deviation and $E_i$ is a deviation for all $i\in\N_n$. If a function $f:D\to I$ satisfies the 
$n$-variable inequality
\Eq{*}{
f\big(\D^G(x_1,\dots,x_n)\big)
\leq\D^E\big(f(x_1),\dots,f(x_n)\big),\qquad(x_1,\dots,x_n\in D),
}
then, for all $k\in\N_n$ and for all injective function $\chi:\N_k\to\N_n$, it also satisfies the 
$k$-variable inequality
\Eq{*}{
f\big(\D^{G_{\chi}}(x_1,\dots,x_k)\big)
\leq\D^{E_{\chi}}\big(f(x_1),\dots,f(x_k)\big),\qquad (x_1,\dots,x_k\in D).
}}

\Rem{*}{Obviously, if, for all $i\in\N_n$, we have $\omega_i=1$ and $E_i(u,v):=u-v$ for all $u,v\in 
I$ in \cor{1K}, or if $X$ is an inner product space, and, for all $i\in\N_n$, we have 
$G_i(x,y)(\cdot):=\ipl x-y,\cdot\ipr$ and 
$E_i(u,v):=u-v$ for all $x,y\in D$ and for all $u,v\in I$, respectively, in \cor{2K}, then, in both 
cases, we get back the reducibility of the Jensen inequality.}

In particular, by applying the previous corollary to the function $f(x)=x$, we immediately
obtain the following consequence for the comparison of deviation means.

\Cor{3K}{Let $I\subseteq\R$ be an interval and $n\in\N$. 
Let further $G,E:I\times I\to\R^n$ such that $G_i$ and $E_i$ are deviations for all 
$i\in\N_n$. If the $n$-variable inequality
\Eq{*}{
\D^G(x_1,\dots,x_n)
\leq\D^E(x_1,\dots,x_n),\qquad(x_1,\dots,x_n\in D)
}
holds, then, for all $k\in\N_n$ and for all injective function $\chi:\N_k\to\N_n$, we also have
the $k$-variable inequality
\Eq{*}{
\D^{G_\chi}(x_1,\dots,x_k)
\leq\D^{E_{\chi}}(x_1,\dots,x_k),\qquad (x_1,\dots,x_k\in D).
}}

The following result establishes the reducibility of an abstract Hölder--Minkowski type inequality.

\Thm{rp1}{
Let $X_1,\dots,X_\ell$ be real Hausdorff topological linear spaces, let $D_1\subseteq 
X_1,\dots,D_\ell\subseteq X_\ell$ be nonempty convex sets and $I\subseteq\R$ be an interval. Let $n\in\N$, 
$k\in\N_n$, and let $\chi:N_k\to\N_n$ be an injective function. Let $N_1:D_1^n\to 
X_1,\dots,N_\ell:D_\ell^n\to X_\ell$ be $\chi$-reducible means and let $M:I^n\to\R$ be a 
$\chi$-continuous, uniquely $\chi$-reducible mean. If a function $f:D_1\times\cdots\times D_\ell\to I$ 
satisfies the $n\cdot\ell$-variable inequality
\Eq{rpn}{
M\big(f(x^1,\dots,x^\ell)\big)\leq f\big(N_1(x^1),\dots,N_\ell(x^\ell)\big),
\qquad (x^1\in D_1^n,\dots,x^\ell\in D_\ell^n),
}
then, for any $\chi$-reductions $K_1:D_1^k\to X_1,\dots,K_\ell:D_\ell^k\to X_\ell$ of 
$N_1,\dots,N_\ell$, respectively, it also fulfills the $k\cdot\ell$-variable inequality
\Eq{rpk}{
M_\chi\big(f(x^1,\dots,x^\ell)\big)\leq f\big(K_1(x^1),\dots,K_\ell(x^\ell)\big),
\qquad (x^1\in D_1^k,\dots,x^\ell\in D_\ell^k),
}
where, for $m\in\N$ and $x^1\in D_1^m,\dots,x^\ell\in D_\ell^m$, we denote
\Eq{*}{
f(x^1,\dots,x^\ell):=(f(x^1_1,\dots,x^\ell_1),\dots,f(x^1_m,\dots,x^\ell_m)).
}}

\begin{proof}
Let $x^1\in D_1^k,\dots,x^\ell\in D_\ell^k$ be arbitrarily fixed, $K_1:D_1^k\to 
X_1,\dots,K_\ell:D_\ell^k\to X_\ell$ be any $\chi$-reduction of $N_1,\dots,N_\ell$, 
respectively, denote $u_1:=K_1(x^1),\dots,u_\ell:=K_\ell(x^\ell)$, finally let 
$u:=(u_1,\dots,u_\ell)$. Using inequality \eq{rpn}, we get 
\Eq{*}{
M\big((f(x^1,\dots,x^\ell)|\chi)(f(u))\big)
&\leq f\big(N_1((x^1|\chi)(u_1)),\dots,N_\ell((x^\ell|\chi)(u_\ell))\big) \\[1mm]
&=f\big(K_1(x^1),\dots,K_\ell(x^\ell)\big)=f(u),
}
that is, the inequality
\Eq{*}{
m_{f(x^1,\dots,x^\ell),M}(f(u))-f(u)=M\big((f(x^1,\dots,x^\ell)|\chi)(f(u))\big)-f(u)\leq 0
}
holds. The mean $M$ is $\chi$-continuous and uniquely $\chi$-reducible, thus, using \lem{4} for 
the vector $x:=f(x^1,\dots,x^\ell)$ and for $y:=f(u)$, we obtain that
\Eq{*}{
M_\chi\big(f(x^1,\dots,x^\ell)\big)\leq 
f(u)=f\big(K_1(x^1),\dots,K_\ell(x^\ell)\big),
}
which finishes the proof.
\end{proof}

To derive various consequences of \thm{rp1}, one can specialize the means $M$ and $N_1,\dots,N_\ell$ 
by letting them equal to a weighted arithmetic mean or to a generalized deviation mean. Then the two choices 
$f(x^1,\dots,x^\ell):=x^1+\cdots+x^\ell$ and $f(x^1,\dots,x^\ell):=x^1\cdots x^\ell$ yield inequalities of 
Minkowski and of Hölder type, respectively.


\end{document}